\documentclass[10 pt]{article}
\usepackage{graphicx}
\usepackage{color}
\usepackage{float,latexsym,amssymb,amsmath}
\usepackage{latexsym,amssymb,amsmath}
\usepackage[T1]{fontenc}

\newenvironment{TH4.1}{\vspace{6pt}\noindent{\bf Proof
of Theorem 4.1:}} { \rule{0.075in}{0.075in} \vspace{11pt}}
\newenvironment{TH4.2}{\vspace{6pt}\noindent{\bf Proof
of Theorem 4.2:}} { \rule{0.075in}{0.075in} \vspace{11pt}}
 \title{\bf   Mixed solutions of monotone iterative technique for hybrid fractional differential equations  }

\author{{\bf  $^{1}$Rabha W. Ibrahim, $^{2*}$Adem Kili\c{c}man, $^{2}$Faten H. Damag }
\\\small  $^{1}$Institute of Mathematical Sciences, University Malaya, 50603, Malaysia
\\\small $^{2}$Department of Mathematics,
University Putra Malaysia, Serdang, Malaysia,
\\\small E-mail: rabhaibrahim@yahoo.com,
akilicman@science.upm.edu.my, faten$_{-}$212326@hotmail.com}
\date{}

\begin{document}
\def\themyth{\thesection.\arabic{myth}}
\def\therem{\thesection.\arabic{rem}}
\def\theexa{\thesection.\arabic{exa}}
\def\thedefn{\thesection.\arabic{defn}}
\def\psfig{\thesection.\arabic{figure}}
\maketitle

{\small {\bf Abstract:} This paper concerns with a mathematical modelling
of biological experiments, and its influence on our lives. Fractional hybrid iterative differential
 equations are equations that interested in  mathematical model of biology. Our technique
 is based on the Dhage fixed point theorem. This tool describes
  mixed solutions by monotone iterative technique in the nonlinear analysis. This method is used to
  combine two solutions: lower and upper. It is shown an approximate result for the hybrid fractional
    differential equations iterative in the closed assembly formed by the lower and upper solutions.

 \vspace{.2cm} {\bf
Keywords:} Fractional differential equation; fractional differential
operator; fractional calculus; monotonous sequences; mixed solutions

\vspace{.2cm}
 }

\section{ Introduction}
Calculus of fractional order power is a field of mathematical
analysis (nonlinear part). It follows the traditional definition of
derivatives and integrals of calculation operators in form
fractional
order(\cite{loverro2004fractional},\cite{podlubny1998fractional},\cite{millerintroduction}).
 Using fractional order differential operator in mathematics modeling has become more and more extended in the last years.
fractional order differential equations have been the concentrate of several studies because of their common occurrence
 in diverse applications in economics, biology, physics and engineering. Recently, a wealth of
 literature developed on the applying nonlinear differential equations of fractional order \cite{kulish2002application}.

\bigskip \noindent The class of fractional order differential equations is a generalization of the class of ordinary
 differential equations. We argue
  that the  fractional order differential equations are more appropriate than the ordinary in mathematical modeling of biological,
   economic and social systems \cite{havlin1995fractals}.
   Fractional calculus is utilized in biology and medicine to explore the potential of fractional differential
  equations to describe and understand the biological organisms grow. Moreover, it utilized to develop the
  structure and functional properties of populations. Extend this concept to evaluate the changes associated with
   the disease hope that contribute to the understanding of the pathogenic processes of medicine \cite{nonnenmacher2013fractals}.
Humans have learned how to employ bacteria and other microbes to
making something useful, such as genetically engineered human
insulin \cite{smith2007bacterial}.

\bigskip \noindent The important of the differential equations of the type hybrid implies polls
number of dynamical systems dealt as special cases,
(\cite{dhage2012basic}, \cite{lu2013theory}). Dhage, Lakshmikantham and Jadhav proved some of the major
 outcomes of hybrid linear differential equations of the first order and second type disturbances
 (\cite{dhage2014approximation},\cite{dhage2013basic},\cite{Dhage2010414}). A great a mathematical
 model for bacteria from growing by the iterative difference equation described. Ibrahim \cite{ibrahim2012}
 established of the existence of an iterative fractional differential equation (Cauchy type)
 using the technique of nonexpansive operator. This kind is created in \cite{ibrahim2015existence}.

\bigskip \noindent In this work, we discuss a mathematical model
of biological experiments, and how its influence on our lives. The
most prominent influence of biological organisms that
 is affect negative or positive in our lives like a bacteria. Fractional hybrid iterative differential
 equations are equations that interested in  mathematical model of biology. Our technique
 is based on the Dhage fixed point theorem. This tool describes
  mixed solutions by monotone iterative technique in the nonlinear analysis. This method is used to
  combine two solutions: lower and upper. It is shown an approximate result for the hybrid fractional
    differential equations iterative in the closed assembly formed by the lower and upper solutions.

\section{Preliminaries}
Recall the following preliminaries:

\bigskip \noindent \textbf{Definition 2.1}
The derivative of fractional ($\gamma$) order for the function
$\phi(s)$ where $0<\gamma<1$ is introduced by
\begin{equation}
D_{a}^{\gamma}\phi(s)=\frac{d}{ds}\int_{a}^{s} \frac{(s-\beta)^
{-\gamma}}{\Gamma(s-\beta)}\phi(\beta)d\beta = \frac{d}{ds} I^{1-\gamma}_{a}\phi(s)\label{111}
\end{equation}
$$(\kappa-1)< \gamma < \kappa,$$
in which $\kappa$ is a whole number and $\gamma$ is real number.

\bigskip \noindent \textbf{Definition 2.2}
The integral of fractional ($\gamma$) order for the function $\phi(s)$ where $\gamma > 0$ is introduced  by
\begin{equation}\label{inegral}
I_{a}^{\gamma}\phi(s)=\int_{a}^{s}\frac{(s-\beta)^{\gamma-1}}{\Gamma(\gamma)}\phi(\beta)d\beta
\end{equation}
While $a = 0$, it becomes $I_{a}^{\gamma}\phi(s)= \phi(s)* \Upsilon_{\gamma}(s)$, wherever $(*)$ signify the convolution product
$$\Upsilon_{\gamma}(s)=\frac{s^{\gamma-1}}{\Gamma(\gamma)}$$
and $\Upsilon_{\gamma}(s)= 0,\,\,\, s \leq 0$ and $\Upsilon_{\gamma}\rightarrow \delta(s) as \gamma \rightarrow 0$
wherever $\delta(s)$ is the delta function

\bigskip \noindent Based on the Riemann-Liouville differential
operator, we impose the following useful definitions:

\bigskip \noindent \textbf{Definition 2.3}
Assume  the closed period bounded interval $I = [s_{0}, s_{0} + a]$
in $\Re$ ($\Re$ the real line), for some $s_{0} \in \Re$, $a\in
\Re.$ The  problem of initial value  of  fractional iterative hybrid
differential equations ( $FIHDE$) can be formulated as

\begin{equation}\label{FIHDE}
D^{\alpha}[v(s)-\psi(s,v(s),v(v(s))]=\aleph(s,v(s),v(v(s))),   s\in
I
\end{equation}
with $v(s_{0})= v_{0}$, where $\psi, \aleph :I \times \Re\rightarrow
\Re$ are continuous. A solution  $v \in C(I ,\Re)$ of the $FIHDE$
\eqref{FIHDE} can be problem by

\begin{enumerate} \item  $s \rightarrow v-\psi(s,v,v(v)))$ is a function which is continuous $\forall v \in \Re$, and
\item $v$ contented the equations in \eqref{FIHDE}.
In which $C(I ,\Re)$ space is of real-valued continuous functions
defined on $I$ .\end{enumerate}

The definitions of the lower and upper solutions of \eqref{FIHDE} as
follows:
\\

\bigskip \noindent \textbf{Definition 2.4}
We said that $\imath\in C(I, \Re)$ is a function which is a lower
solution for the equation introduced on $I$ if

\begin{enumerate} \item $ s \mapsto (\imath(s )-\psi(s, \imath(s)), \imath(\imath(s))))$, is continuous, and
 \item  $D^{\alpha}[\imath(s) -\psi(s, v(s),v(v(s)))]\geq \aleph(s,\imath(s),\imath(\imath(s))), s \in I ,\imath(s_{0})\geq v_{0}.$
 \end{enumerate}

\bigskip \noindent \textbf{Definition 2.5}
We said that  $\tau\in C(I, \mathbb{R})$ is a function which is an
upper solution for the equation introduced on $I$ if

\begin{enumerate}
\item  $ s \mapsto ( \tau(s)-\psi(s, \tau(s), \tau(\tau(s)))$, is continuous, and
\item $D^{\alpha}[\tau(s)-\psi(s,v(s),v(v(s)))]\leq \aleph(s,\tau(s),\tau(\tau(s))), s\in I,\tau(s_{0})\leq v_{0}.$
\end{enumerate}

\bigskip \noindent We can build the monotonous sequence of consecutive iterations to converging
towards the extremes among the lower and upper solutions of the
differential equation related hybrid on$I$. We treat the case that
if $\psi$ is neither non-decreasing nor non-increasing in the state
of the variable $v$. If the function $\aleph$ can be separated into
two components

$$\aleph(s,v,v(v)))=\aleph_{1}(s,v,v(v)))+\aleph_{2}(s,v,v(v))$$
where $\aleph_{1}(s,v,v(v)))$ is a non-decreasing component while
 another component is not $\aleph_{2}(s,v,v(v)))$ increases in the state variables of $v$,
  then we may be constructed sequences iteration converged to solutions extremal $FIHDE$\eqref{FIHDE} on $I$.\\

\bigskip \noindent \textbf{Definition 2.6} Currently thought to be a initial value problem $FIHDE$
\begin{equation}
\left\{ \begin{array}{c}
     D^{\alpha}[v(s)-\psi(s,v(s),v(v(s))]=\aleph_{1}(s,v,v(v)))+\aleph_{2}(s,v,v(v))), s\in
I,\\
        v(s_{0})=v_{0}
        \end{array}\right.\label{hde}
\end{equation}
where, $\psi \in C(I \times R, R)$ and $\aleph_{1}, \aleph_{2} \in \mathfrak{L}(I \times R, R)$.\\

Thus the lower and upper solutions of \eqref{hde} can be as defined
as follows:

\bigskip \noindent \textbf{Definition 2.7} The functions $\sigma,\rho \in C(I ,\Re)$ fulfill
the following condition: the maps
 $s \rightarrow \sigma(s )- \psi(s
,\sigma(s),\sigma(\sigma(s)))$ and $s \rightarrow \rho(s )-
\psi(s ,\rho(s),\rho(\rho(s)))$ are absolute continuous on $I$. Thus the functions $(\sigma,\rho)$ are supposed to be of the kind\\

\noindent  (a) which is mixed lower solutions and upper solutions
for  \eqref{hde} on $I$, sa following
\begin{equation}
\left\{ \begin{array}{c}
     D^{\alpha}[\sigma(s)-\psi(s,\sigma(s),\sigma(\sigma(s))]\leq\aleph_{1}(s, \sigma, \sigma(\sigma(s)))+\aleph_{2}(s,\rho(s),\rho(\rho(s)))),   s\in
I,\\
        \sigma(s_{0})\leq v_{0}
        \end{array}\right.\label{hde1}
\end{equation}
and
\begin{equation}
\left\{ \begin{array}{c}
     D^{\alpha}[\rho(s)-\psi(s,\rho(s),\rho(\rho(s))]\geq \aleph_{1}(s, \rho, \rho(\rho(s)))+\aleph_{2}(s,\sigma(s),\sigma(\sigma(s))),   s\in
I,\\
        \rho(s_{0})\geq v_{0}
        \end{array}\right.\label{hde2}
\end{equation}
  Whether the sign was of equality achieves in relationships \eqref{hde1} and \eqref{hde2},hence the even of functions $(\sigma,\rho)$ set is been calling a mixed solution of kind (a) for the $FIHDE$ \eqref{hde} on $I$. \\

\noindent  (b) which is mixed lower solutions and upper for
\eqref{hde} on $I$, as follows
\begin{equation}
\left\{ \begin{array}{c}
     D^{\alpha}[\sigma(s)-\psi(s,\sigma(s),\sigma(\sigma(s))]\leq\aleph_{1}(s, \rho, \rho(\rho(s)))+\aleph_{2}(s,\sigma(s),\sigma(\sigma(s)))),   s\in
I,\\
        \sigma(s_{0})\leq v_{0}
        \end{array}\right.\label{hde3}
\end{equation}
and
\begin{equation}
\left\{ \begin{array}{c}
     D^{\alpha}[\rho(s)-\psi(s,\rho(s),\rho(\rho(s))]\geq \aleph_{1}(s, \sigma, \sigma(\sigma(s)))+\aleph_{2}(s,\rho(s),\rho(\rho(s))),   s\in
I,\\
        \rho(s_{0})\geq v_{0}
        \end{array}\right.\label{hde4}
      \end{equation}
  Whether the sign was of equality achieves in relationships \eqref{hde3} and \eqref{hde4}, hence
   the even of functions $(\sigma, \rho)$ set is been calling a mixed solution of kind (b) for the \eqref{hde} on $I$.\\


\subsection{Assumptions}
In the following assumptions relating to function $\psi$ is  very important in the studying of Eq\eqref{hde}.\\

(a0) The function $v \mapsto (v - \psi(s_{0}, v, v(v)))$ is injective in $\Re$.\\

(b0) $\aleph$ is a bounded real-valued function on $I \times\Re$.\\

(a1)\,  The function $v \mapsto (v-\psi(s, v, v(v)))$ is increasing in $\Re$ for all $s\in I$ .\\

(a2) \, There is a constant $\ell >0$ so that \[| \psi (s, v, v(v))-
\psi (s, z, z(z))|\leq \frac{\ell|v - z|}{M+|v - z|}, \quad M > 0,\]
 $\forall s \in I,$\,  $v,z\in \Re$ and $\ell \leq M$.\\

(b1)\,  There is a constant $\kappa >0$ so that $|\aleph(s,v,v(v)|
\leq \kappa$  $\forall s \in I$ and  $\forall v \in \Re$.\\

(b2) \, $\aleph_{1}(s, v, v(v))$ is function which is non-decreasing
in $v$ function, and $\aleph_{2}(s, v, v(v))$
 is function which is not increasing in
$v$ for each $s\in I$\\

(b3)\,   $(\sigma_{0}, \rho_{0})$ is Functions which are mixing the
lower and upper solutions for  \eqref{hde}
 kind(a) on $I$ with $\sigma_{0}\leq \rho_{0}$
.\\

(b4) \, The pair is  $(\sigma_{0}, \rho_{0})$, the upper and lower
mixing solutions for  \eqref{hde} kinds (b) on $I$
with $\sigma_{0}\leq \rho_{0}.$\\


\section{Main results} \bigskip \noindent

In this section, our purpose is to discuss the approximation outcome
for  \eqref{hde}.\\

\bigskip \noindent \textbf{Lemma 3.1}(\cite{lu2013theory})
Suppose the assumptions $(a0) - (b0)$ are achieved. Then the
function $v$ is a solution for Eq.\eqref{FIHDE} if and only if it
must be the solution of the fractional iterative of hybrid equation
integrated $FIHIE$
\begin{equation}\label{FIHIE}
v(t)=[v_{0}-\psi(s_{0},v_{0},v(v_{0}))]+\psi(s,v(s),v(v(s)))+\int_{0}^{s}
\aleph(\beta,v(\beta),v(v(\beta)))\frac{(s-\beta)^{\alpha-1}}{\Gamma(\alpha)}d\beta,
\end{equation}
\[(s\in I, \,\, v(0)=v_0).\]

\bigskip \noindent \textbf{Theorem 3.1}\,
(\cite{dhage2004fixed})\, Let $\varrho$ be a closed convex and
bounded subset of the Banach space $A.$ Moreover, let
$Q : A \rightarrow A$ and $P : \varrho \rightarrow A$ be two operators so that\\

(i) $Q$ is nonlinear D-contraction,\\

(ii) $P$ is compact and continuous, \\

(iii) $v = Qv +Pz$ for all $v \in \varrho \Rightarrow z \in
\varrho$.

\bigskip \noindent \textbf{Theorem3.2} \, \label{re1} Let the assumptions $(a1),(a2)$ and $(b1$) be hold.
Then \eqref{FIHDE} has a solution on $I$.

\bigskip \noindent \textbf{Proof.}
Let $A = C(I ,\Re)$ be a set and  $\c{c} \subseteq A,$ such that

\begin{equation}
\varrho = \{v\in A |\| A\|\leq M\}
\end{equation}
where,
$$M =|v_{0} - \psi(s_{0}, v_{0}, v(v(0))|+\ell+\Psi_{0} +\frac{a^{\alpha}}{\Gamma(\alpha+1)} \|\xi\|_{\ell^{1}}.$$
and $\Psi_{0} = \sup_{s\in I} \mid \psi (s ,0,0)\mid.$ Obviously
$\varrho$ is a convex, bounded and closed subset of the  space $A$.
By using  the assumptions (a1) and (b1) together with the help of
the Lemma 3.1, we conclude that the $FIHDE$\eqref{FIHDE} is
tantamount to the nonlinear $FIHIE$\eqref{FIHIE}.  We  define two
operators $Q : A \rightarrow A$   and $ P : \varrho \rightarrow A$
as follows:

\begin{equation}
Qy(s )= \psi (s , v(s ),v(v(s))), s\in I,
\end{equation}
and
\begin{equation}\label{23}
Pv(s)=[v_{0}-\psi(s_{0},v_{0},v(v_{0}))]+\int_{0}^{s}\aleph(\beta,v(\beta),v(v(\beta)))\frac{(s-\beta)^{\alpha-1}}{\Gamma(\alpha)}d\beta,
s\in I.
\end{equation}

Consequently, the $FIHIE$\eqref{FIHIE} is equivalent to the operator
equation

\begin{equation}\label{24}
Qv(s )+Pv(s )= v(s ), s\in I.
\end{equation}
We demonstrate that the operators $Q$ and $P$ fulfill all the
conditions of Theorem 3.1. Foremost, we examine that $Q$ is a
nonlinear $\Upsilon$-contraction on $Q$ with a $\Upsilon$ function
$\varphi$. Let $v, z \in A$. In view of assumption $(a2)$, we
conclude that
$$|Qv(s )- Qz(s )| =| \psi (s , v(s ))- \psi (s , z(s ))|\leq \frac{\ell|v(s )- z(s)|}{M+|v(s) - z(s)|}\leq \frac{\ell|v - z|}{M+|v - z|}$$
for all $s \in I$. Take the supremum over $s$ yields

$$\| Av - Az\|\leq \frac{\ell|v - z|}{M+|v - z|}$$
 $\forall v, z \in A$. This proves that $Q$ is a nonlinear $D$-contraction $A$ with the $D$-function $\varphi$ defined
by $\varphi(r ) =\frac{\ell r}{M+r}$.

\bigskip \noindent Next, we examine that $P$ is a continuous  and compact
operator on $\varrho$ into $A$. Let $\{v_{t}\}$ be a sequence in
$\varrho$ converging to a point $v \in \varrho,$ thus we have

$$\lim_{t\rightarrow \infty} Pv_{t}(s)=\lim_{t\rightarrow \infty}[v_{0}-\psi(s_{0},v_{0},v(v_{0}))
+\int_{0}^{s}\aleph(\beta,v_{t}(\beta),v_{t}(v_{t}(\beta)))\frac{(s-\beta)^{\alpha-1}}{\Gamma(\alpha)}
d\beta]$$

$$=v_{0}-\psi(s_{0},y_{0},v(v_{0}))+\lim_{t\rightarrow \infty} \int_{0}^{s}\aleph(\beta,v_{t}
(\beta),v_{t}(v_{t}(\beta)))\frac{(s-\beta)^{\alpha-1}}{\Gamma(\alpha)}d\beta$$

$$ =
 v_{0}-\psi(s_{0},v_{0},v(v_{0}))+ \int_{0}^{s}\lim_{t\rightarrow \infty}[\aleph(\beta,v_{t}
 (\beta),v_{t}(v_{t}(\beta)))\frac{(s-\beta)^{\alpha-1}}{\Gamma(\alpha)}]d\beta$$

 $$= v_{0}-\psi(s_{0},v_{0},v(v_{0}))+ \int_{0}^{s}\aleph(\beta,v(\beta),v(v(\beta)))\frac{(s-\beta)^{\alpha-1}}{\Gamma(\alpha)}d\beta=Pv(s)$$
for all $s \in I$. Now, we proceed to prove that $\{Pv_{t}\}$ is
equi-continuous with respect to $v$. According to \cite{Granas1991R},
we attain that $P$ is a continuous  operator on $\varrho$. To show
that $P$ is a compact operator on $\varrho$. It suffices to examine
that $\varrho$ is a regularly bounded and equi-continuous set in $A$.
Let $v\in \varrho$ be arbitrary, then by the assumption (b1), we
have

$$|Pv(s)|\leq|v_{0}-\psi(s_{0},v_{0},v(v_{0}))|+\int_{0}^{s}|\aleph(\beta,v(\beta),v(v(\beta)))\frac{(s-\beta)^{\alpha-1}}{\Gamma(\alpha)}|d\beta$$

$$\leq |v_{0}-\psi(s_{0},v_{0},v(v_{0}))|+\int_{0}^{s}\xi(\beta)
\frac{(s-\beta)^{\alpha-1}}{\Gamma(\alpha)}d\beta \leq
|v_{0}-\psi(s_{0},v_{0},v(v_{0}))|+\frac{a^{\alpha}}{\Gamma(\alpha+1)}\|\xi\|_{\ell^{1}}$$
for all $s\in I$ . By taking the supremum over $t$, we obtain

$$|Pv(s)|\leq|v_{0}-\psi(s_{0},v_{0},v(v_{0}))|+ \frac{a^{\alpha}}{\Gamma(\alpha+1)}\|\xi\|_{\ell^{1}}$$
$ \forall v \in \varrho$. This proves that $P$ is uniformly bounded
on $\varrho$.\\ Also let  $s_{1}, s_{2}\in I $ with $s_{1}<s_{2}$.
Then for any $v \in \varrho$, one has
$$|Pv(s_{1})- Pv(s_{2})|= |\int_{{s_{0}}}^{s_{1}}|\aleph(\beta,v(\beta),v(v(\beta)))\frac{(s_{1}-\beta)^{\alpha-1}}
{\Gamma(\alpha)}d\beta-\int_{s_{0}}^{s_{2}}|\aleph(\beta,v(\beta),v(v(\beta)))\frac{(s_{2}-\beta)^{\alpha-1}}{\Gamma(\alpha)}d\beta
|$$
$$\leq |\int_{{s_{0}}}^{s_{1}}|\aleph(\beta,v(\beta),v(v(\beta)))\frac{(s_{1}-\beta)^{\alpha-1}}{\Gamma(\alpha)}
d\beta-\int_{s_{0}}^{s_{1}}|\aleph(\beta,v(\beta),v(v(\beta)))\frac{(s_{2}-\beta)^{\alpha-1}}{\Gamma(\alpha)}d\beta
|$$
$$+ |\int_{{s_{o}}}^{s_{1}}|\aleph(\beta,v(\beta),v(v(\beta)))\frac{(s_{2}-\beta)^{\alpha-1}}{\Gamma(\alpha)}d\beta
-\int_{s_{0}}^{s_{2}}|\aleph(\beta,v(\beta),v(v(\beta)))\frac{(s_{2}-\beta)^{\alpha-1}}{\Gamma(\alpha)}d\beta
|$$
$$\leq \frac{\|\xi\|_{\ell^{1}}}{\Gamma(\alpha+1)}[|(s_{2}-s_{2})^{\alpha}-(s_{1}-s_{0})^{\alpha}-(s_{2}-s_{1})^{\alpha}|+(s_{2}-s_{1})^{\alpha}]$$
Hence, for $\delta > 0$, there exists a $\epsilon > 0$ so that

$$|s_{1}-s_{2}|< \epsilon \Rightarrow |Pv(s_{1})-Pv(s_{2})| < \delta$$
 $\forall s_{1}, s_{2} \in I$ and  $ \forall v\in
\varrho$. This  examines for $P(\varrho)$ is equi-continuous in $A$.
presently  $P(\varrho)$ is bounded and hence it is compact by
Arzel$\`{a}$-Ascoli Theorem. Resulting, $\varrho$ is a continuous and compact operator on $\varrho$.\\
Then, we prove that assumptions (iii) of Theorem 3.1 is fulfilled.
Let $v \in A$ be fixed and $z \in \varrho$ be arbitrary such that $v
= Qv +Pz$. In view of the assumption (a2) yields

$$|v(s )|\leq|Qv(s )|+|Pz(s )| $$

$$\leq |v_{0}- \psi (s_{0}, v_{0})|+| \psi (s , v(s ),v(v(s))|+\int_{0}^{s}|\aleph(\beta,v(\beta),v(v(\beta)))
\frac{(s-\beta)^{\alpha-1}}{\Gamma(\alpha)}|d\beta$$

$$\leq |v_{0} - \psi (s_{0}, v_{0})|+| \psi (s , v(s ),v(v(s))|+\int_{0}^{s}|\aleph(\beta,v(\beta),v(v(\beta)))\frac{(s-\beta)^{\alpha-1}}
{\Gamma(\alpha)}|d\beta$$

$$ \leq |v_{0}- \psi (s_{0}, v_{0}, v(v_{0}))|+\ell+\Psi_{0}+\int_{0}^{s}|
\xi(\beta)\frac{(s-\beta)^{\alpha-1}}{\Gamma(\alpha)}|d\beta$$

$$\leq
|v_{0}-\psi(s_{0},v_{0}, v(v_{0}))|+\ell+\Psi_{0}+
\frac{a^{\alpha}}{\Gamma(\alpha+1)}\|\xi\|_{\ell^{1}}.$$

Take the supremum over $s$, implies
$$\|v\|\leq|v_{0}-\psi (s_{0}, v_{0}, v(v_{0}) )|+\ell+\Psi_{0} +
\frac{a^{\alpha}}{\Gamma(\alpha+1)}\|\xi\|_{\ell^{1}}=M.$$
Thus, $v \in \varrho.$\\
Therefore, fulfilled all conditions of the Theorem 3.1 and thus the
operator equation $v = Qv +Pz$ has a solution in $\varrho$.
Resulting, the $FIHDE$\eqref{FIHDE} has a solution introduced  on
$I$. This completes the proof.


\bigskip \noindent \textbf{Theorem 3.3} Let $\imath, \tau \in C(I ,\Re)$ be lower and upper solutions of
$FIHDE$\eqref{FIHDE} fulfilling $\imath(s )\leq \tau(s ),s \in I $
and let the assumptions $(a1)-(a2)$ and $(b1)$ achieved. Then, there
is a solution $v(s )$ of \eqref{FIHDE}, in the closed set  $
\overline{\mho}$, satisfying \[\imath(s )\leq v(s )\leq \tau(s ), s
\in I.
\]

\bigskip \noindent \textbf{Proof.}
Assume that $\Theta : I \times \Re \mapsto \Re $ is a function
defined by

$$\Theta(s , v, v(v))= \max\{\imath(s ),\min{v(s ), \tau(s )}\},$$
satisfying \[\check{ \aleph} (s ,v , v(v))):= \aleph(s ,\Theta(s , v
,v(v)))).\]
 Moreover, define a continuous extension of $\aleph $ on
$I\times\Re$ such that
 \[\check{| \aleph} (s, v, v(v)))|= |\aleph(u
,\Theta(s , v,v(v))))| \leq \kappa , s \in I\,\, \forall v \in
\Re.\]
  In view of Theorem 3.2, the $FIHDE$

\begin{equation}
\left\{ \begin{array}{c}
D^{\alpha}[v(s)-\psi(s,v(s),v(v(s))] =\check{\aleph} (s ,v,v(v))), s\in
I\\\\
v(u_{0}) = v_{0}\in \Re
\end{array}\right.\label{1}
\end{equation}
has a solution $v$ defined on $I$ .\\
For any $\delta >0 $, define

\begin{equation}
\imath_{\delta}(s )− \psi (s ,\imath_{\delta}(\imath_{\delta}(s
)))= (\imath(s )-
\psi(s,\imath(s),\imath(\imath(s)))−\delta(1+s)\label{2}
\end{equation}
and

\begin{equation}
\tau_{\delta}(s )− \psi (s ,\tau_{\delta}(\tau_{\delta}(s )))=
(\tau(s)- \psi(s,\tau(s),\tau(\tau(s)))−\delta(1+s)\label{3}
\end{equation}for $s\in I$ . In virtue of the assumptions $(a1)$, we get

\begin{equation}
\imath_{\delta}(s )< \imath(s),\,\,\ and,\,\,\ \tau(s)<
\tau_{\delta}(s )\label{4}
\end{equation}
for $s \in I.$ Since

$$\imath(s_{0})\leq v_{0} \leq \tau(s_{0}),$$
one has
\begin{equation}
\imath_{\delta}(s_{0} )< v_{0} < \tau_{\delta}(s_{0}).\label{5}
\end{equation}
To show that
\begin{equation}
\imath_{\delta}(s )< v_{0} < \tau_{\delta}(s),  \,\ s\in I
\label{6},
\end{equation}
we define \[v(s )= v(s )- \psi(s,v(s),v(v(s)), s \in I. \] Likewise,
we consider

$$\hbar_{\delta}(s ) =\imath_{\delta}(s)- \psi(s ,\imath_{\delta}(\imath_{\delta}(s))),$$
 $$ \hbar(s ) =\imath(u)- \psi(s,\imath(s),\imath(\imath(s)),$$
and
$$T_{\delta}(s)=\tau_{\delta}(s) \psi(s,\tau_{\delta}(s),\tau(\tau_{\delta}(s)),$$
$$ T(s)= \tau(s) \psi(s,\tau(s),\tau(\tau(s))$$
 $\forall s \in I $.
If Eq.\eqref{6} is wrong, then there  exists  a $s_{\varepsilon} \in
(s_{0},s_{0}+a]$ such that
$$v(_{\varepsilon} ) = \tau_{\delta}(s_{\varepsilon} )$$
and

$$\imath_{\delta}(s ) < v(s ) < \tau_{\delta}(s ), s_{0} \leq s < s_{\varepsilon}$$
If $v(s_{\varepsilon}) > \tau(s_{\varepsilon})$, then
$\Theta(s_{\varepsilon}, v(s_{\varepsilon}),v(v(s_{\varepsilon}))) =
\tau(s_{\varepsilon})$. Furthermore,
$$\imath(s_{\varepsilon}) \leq \Theta(s_{\varepsilon}, v(s_{\varepsilon}),v(v(s_{\varepsilon}))) \leq \tau(s_{\varepsilon}).$$
Now,
$$ D^{\alpha}T(s_{\varepsilon}) \geq \aleph(s_{\varepsilon}, \tau(s_{\varepsilon}),\tau (\tau(s_{\varepsilon}))) = \check{ \aleph}
 (s_{\varepsilon} ,v(s_{\varepsilon}),v(v(s_{\varepsilon}))))= D^{\alpha}V(s )$$
 $ \forall s \in I $. Since
$T_{\delta}(us )>D^{\alpha}T(s )$,
 $ \forall s\in I$ , we have

\begin{equation}
D^{\alpha}T_{\delta} (s_{\varepsilon})>
D^{\alpha}V(s_{\varepsilon}).\label{7}
\end{equation}
But,
$$V(s_{\varepsilon}) =T_{\delta}(s_{\varepsilon})$$
also
$$V(s )= T_{\delta}(s ), s_{0} \leq s < s_{\varepsilon},$$
means that together
$$\frac{V(s_{\varepsilon}+\rho)-V(s_{\varepsilon})}{\rho^{\alpha}}>
\frac{T_{\delta}(s_{\varepsilon}+\rho)-T_{\delta}(s_{\varepsilon})}{\rho^{\alpha}}$$
if $\rho < 0$ a small. Take the limit $\rho\rightarrow0$ in the up
variance yields

$$D^{\alpha}V(s_{\varepsilon}) \geq D^{\alpha}T_{\delta}(s_{\varepsilon})$$
that is a contradiction to \eqref{7}. Hence,
$$v(s )<\tau_{\delta}(s )$$ $ \forall s \in I $. Consequently
$$\imath_{\delta}(s )< v(s )< \tau_{\delta}(s ), s\in I .$$
Letting $\delta \rightarrow 0$ in the up inequality, we get
$$\imath(s )\leq v(s )\leq \tau(s ), s \in I.$$
This completes the proof.
\,\,\,\,\,\,\,\,\,\,\,\,\,\,\,\,\,\,\,\,\,\,\,\,\,\,\,\,\,\,\,\,\,\,\,\,\,\,$\square$

\bigskip \noindent \textbf{Theorem 3.4}\label{1b}
Let assumptions (a1) - (a2) and (b2) - (b3) achieved. Then there are the monotonous sequences
$\{\sigma_{t}\}, \{\rho{α}_{t}\}$ such that $\sigma_{t}\rightarrow\sigma$ and $\rho_{t} \rightarrow \rho$
uniformly on $I$ in which $(\sigma, \rho)$ are mixed extremal solutions $FIHDE$\eqref{hde} type(a) on $I$.

\bigskip \noindent \textbf{Proof.}
Note the following a quadratic $FIHDE$

\begin{equation}
\left\{ \begin{array}{c}
     D^{\alpha}[\sigma_{t+1}(s)-\psi(s,\sigma_{t+1}(s),\sigma(\sigma_{t+1}(s))]\leq\aleph_{1}(s, \sigma_{t}(s), \sigma(\sigma_{t}(s)))+\aleph_{2}(s,\rho_{t}(s),\rho(\rho_{t}(s)))),   s\in
I,\\
        \sigma_{t+1}(s_{0})\leq v_{0}
        \end{array}\right.\label{hde5}
\end{equation}
and
\begin{equation}
\left\{ \begin{array}{c}
     D^{\alpha}[\rho_{t+1}(s)-\psi(s,\rho_{t+1}(s),\rho(\rho_{t+1}(s))]\geq \aleph_{1}(s, \rho_{t}(s), \rho(\rho_{t}(s)))+\aleph_{2}(s,\sigma_{t}(s),\sigma(\sigma_{t}(s))),   s\in
I,\\
        \rho_{t+1}(s_{0})\geq v_{0}
        \end{array}\right.\label{hde6}
\end{equation}
for $t \in N$.\\

Obviously, the equations \eqref{hde5} and \eqref{hde6} having unique
solutions $\sigma_{t+1}$ and  $\rho_{t+1}$ on $I$ respectively given
Banach contraction mapping principle. We now want to demonstrate
that
\begin{equation}
\sigma_{0} \leq \sigma_{1} \leq \sigma_{2} \leq \ldots \leq\sigma_{t} \leq \rho_{t} \leq \ldots \leq \rho_{2} \leq \rho_{1} \leq \rho_{0}
\end{equation}
on $I$ for $t =0,1,2, \ldots$  Let $t =0$ and set
$$\Theta(s)- \psi (s ,\Theta(s ),\Theta(\Theta(s)))=(\sigma_{0}(s)- \psi (s ,\sigma_{0}(s),\sigma(\sigma_{0}(s)))- )))−\sigma_{1}(s)- \psi (s ,\sigma_{1}(s),\sigma(\sigma_{1}(s)))$$
for $s \in I$. Next by monotonicity of $\aleph_{1}$ and
$\aleph_{2}$, we get
$$D^{\alpha}[\Theta(s)- \psi (s ,\Theta(s ),\Theta(\Theta(s)))]=D^{\alpha}[(\sigma_{0}(s)− \psi (s ,\sigma_{0}(s),\sigma(\sigma_{0}(s)))]-D^{\alpha}[\sigma_{1}(s)-\psi (s ,\sigma_{1}(s),\sigma(\sigma_{1}(s))))]$$
$$\leq \aleph_{1}(s_{0},\sigma_{0}(s ),\sigma(\sigma_{0}(s)))+\aleph_{2}(s ,\rho_{0}(s ),\rho(\rho_{0}(s )))-\aleph_{1}(s_{0} ,\rho_{0}(s ),\rho(\rho_{0}(s )))+\aleph_{2}(s ,\sigma_{0}(s),\sigma(\sigma_{0}(s)))$$
$$=0$$
$ \forall s\in I$ and $\Theta(s_{0})=0$. This implies that
$$\sigma_{0}(s)- \psi (s ,\sigma_{0}(s),\sigma(\sigma_{0}(s)))\leq \sigma_{1}(s)- \psi (s ,\sigma_{1}(s),\sigma(\sigma_{1}(s))),$$
$ \forall s\in I$. In view of (a1), one can get  $\sigma_{0}(s)\leq
\sigma_{1}(s)$, $\forall s\in I$. Likewise it can be demonstrated
which $\rho_{1}(s ) \leq \rho_{0}(s )$ on $I$. Setting
$$\Theta(s)-\psi (s ,\Theta(s ),\Theta(\Theta(s)))=
(\sigma_{1}(s)-\psi (s
,\sigma_{1}(s),\sigma(\sigma_{1}(s))))-(\rho_{1}(s)- \psi (s
,\rho_{1}(s),\rho(\rho_{1}(s))))$$
 $ \forall s\in I$. By
monotonicity of $\aleph_{1}$ and $\aleph_{2}$, we obtain

$$D^{\alpha}[\Theta(s)-\psi (s ,\Theta(s),\Theta(\Theta(s)))]=D^{\alpha} [\sigma_{1}(s)-\psi (s ,\sigma_{1}(s),\sigma(\sigma_{1}(s))))]-D^{\alpha}[(\rho_{1}(s)− \psi (s ,\rho_{1}(s),\rho(\rho_{1}(s))))]$$
$$\leq \aleph_{1}(s_{0},\sigma_{0}(s ),\sigma(\sigma_{0}(s)))+\aleph_{2}(s ,\rho_{0}(s ),\rho(\rho_{0}(s )))-\aleph_{1}(s_{0} ,\rho_{0}(s ),\rho(\rho_{0}(s )))+\aleph_{2}(s ,\sigma_{0}(s),\sigma(\sigma_{0}(s)))$$
$$\leq 0$$
$ \forall s\in I$ and $\Theta(s_{0})=0$. This leads to
$$\sigma_{1}(s )− \psi (s ,\sigma_{1}(s),\sigma(\sigma_{1}(s))) \leq \rho_{1}(s )- \psi (s ,\rho_{1}(s),\rho(\rho_{1}(s)))$$
$ \forall s\in I$. By (a1), we attain to  \[\sigma_{1}(s)\leq
\rho_{1}(s), \quad \forall s\in I.\] Next,  for  $\jmath\in N$,
yields
$$\sigma_{\jmath−1}\leq \sigma_{\jmath} \leq\rho_{\jmath}\leq \rho_{\jmath-1}$$
and hence
$$\sigma_{\jmath}\leq \sigma_{\jmath+1} \leq\rho_{\jmath+1}\leq \rho_{\jmath}.$$
Setting

$$\Theta(s)- \psi (s ,\Theta(s),\Theta(\Theta(s)))=
 (\sigma_{\jmath}(s )-\psi(s ,\sigma_{\jmath}(s ),\sigma(\sigma_{\jmath}(s ))))-(\sigma_{\jmath+1}(s)-\psi(s ,\sigma_{\jmath+1}(s ),\sigma(\sigma_{\jmath+1}(s ))))$$
Then the humdrum of $\aleph_{1}$ and $\aleph_{2},$ we receive

$$D^{\alpha}[\Theta(s)-\psi (s ,\Theta(s),\Theta(\Theta(s)))]= D^{\alpha}[(\sigma_{\jmath}(s )-\psi(s ,\sigma_{\jmath}(s ),\sigma(\sigma_{\jmath}(s ))))]-D^{\alpha}[(\sigma_{\jmath+1}(s)-\psi(s ,\sigma_{\jmath+1}(s ),\sigma(\sigma_{\jmath+1}(s ))))]$$
$$\leq  \aleph_{1}(s ,\sigma_{\jmath−1},\sigma(\sigma_{\jmath-1}(s ) )+\aleph_{2}(s,\rho_{\jmath-1}, \rho(\rho_{\jmath-1}))-\aleph_{1}(s ,\sigma_{\jmath},\sigma(\sigma_{\jmath}))-\aleph_{2}(s,\rho_{\jmath},\rho(\rho_{\jmath}))$$
$$\leq 0$$
$ \forall s\in I$ and $\Theta(s_{0})=0$. This implies that
$$\sigma_{\jmath}-\psi(s ,\sigma_{\jmath}(s ),\sigma(\sigma_{\jmath}(s )))\leq \sigma_{\jmath+1}(s )-\psi(s ,\sigma_{\jmath+1}(s ),\sigma(\sigma_{\jmath+1}(s )))$$
for every $s\in I$. Since assumption (a1) achieved, we have $\sigma_{\jmath}(s)\leq \sigma_{\jmath+1}(s)$, $\forall s\in I$. Likewise it can be demonstrated which $\rho_{\jmath+1}(s ) \leq \rho_{\jmath}(s )$ on $I$ .
The same way it is assumed that the inequality
$$\sigma_{\jmath-1}\leq \sigma_{\jmath} \leq\rho_{\jmath}\leq \rho_{\jmath-1}$$
achieves on $I$. We are going to demonstrate that
$$\sigma_{\jmath}\leq \sigma_{\jmath+1} \leq\rho_{\jmath+1}\leq \rho_{\jmath}$$
on $I$. Set
$$\Theta(s)- \psi (s ,\Theta(s ),\Theta(\Theta(s )))=(\sigma_{\jmath+1}(s)- \psi(s,\sigma_{\jmath+1}(s),\sigma(\sigma_{\jmath+1}(s))))-(\rho_{\jmath+1}(s)-\psi (s ,\rho_{\jmath+1},\rho(\rho_{\jmath+1})))$$
for $s\in I$. So by monotonicity of $\aleph_{1}$ and $\aleph_{2}$ we get
$$D^{\alpha}[\Theta(s)- \psi (s ,\Theta(s ),\Theta(\Theta(s )))]=D^{\alpha}[(\sigma_{\jmath+1}(s)- \psi(s,\sigma_{\jmath+1}(s),\sigma(\sigma_{\jmath+1}(s))))]-D^{\alpha}[(\rho_{\jmath+1}(s)-\psi (s ,\rho_{\jmath+1},\rho(\rho_{\jmath+1})))]$$
$$\leq \aleph_{1}(s ,\sigma_{\jmath}(s),\sigma(\sigma_{\jmath}(s)))+\aleph_{2}(s ,\rho_{\jmath}(s),\rho(\rho_{\jmath}(s)))-\aleph_{1}(s ,\rho_{\jmath+1},\rho(\rho_{\jmath+1}))-\aleph_{2}(s ,\sigma_{\jmath}(s),\sigma(\sigma_{\jmath}(s)))$$
$$\leq 0$$
for the whole $s\in I$ and $\Theta(s_{ 0})=0$. This means that
$$\sigma_{\jmath+1}(s)- \psi(s,\sigma_{\jmath+1}(s),\sigma(\sigma_{\jmath+1}(s))))\leq \rho_{\jmath+1}- \psi (s ,\rho_{\jmath+1},\rho(\rho_{\jmath+1}))$$

for every $s\in I$. Since assumption (a1) is achieved, we have $\sigma_{\jmath+1}(s)\leq \rho_{\jmath+1}(s)$, $\forall s\in I$.\\

Presently it is readily shown that the sequence $\{\sigma\}$ and $\{\rho\}$ are bounded uniformly and equi-continuous sequences and have therefore converge uniformly on $I$. As are monotonous sequences, $\{\sigma_{t}\}$ and $\{\rho_{t}\}$  converse uniformly monotonous $\sigma$ and $\rho$ on $I$ respectively.
Course, the pair $(\sigma,\rho)$ is a mixed solution of these equations \eqref{hde} on $I$. Lastly, we establish which $(\sigma,\rho)$ is a mixed solution of minimum and maximum for the equations \eqref{hde} on $I$. Let $v$ whatever solution of the equations \eqref{hde} on $I$ as $\sigma_{0}(s) \leq v (s) \leq \rho (s)$ on$I$. Assume that for $\jmath \in N$, $\sigma_{\jmath}(s) \leq v (s) \leq \rho_{\jmath} (s)$, $s \in I$. We will demonstrate which $\sigma_{\jmath+1}(s) \leq v (s) \leq \rho_{\jmath+1} (s)$, $s \in I$. adjustment
$$\Theta(s )- \psi (s ,\Theta(s),\Theta(\Theta(s)))= (\sigma_{\jmath+1}(s)- \psi (s ,\sigma_{\jmath+1}(s),\sigma(\sigma_{\jmath+1}(s))))-(v(s )- \psi (s , v(s ),v(v(s))))$$
for every $s\in I$. After, for the monotony of $\aleph_{1}$ and $\aleph_{2}$ we get
$$D^{\alpha}[\Theta(s )- \psi (s ,\Theta(s),\Theta(\Theta(s)))]= D^{\alpha}[(\sigma_{\jmath+1}(s)- \psi (s ,\sigma_{\jmath+1}(s),\sigma(\sigma_{\jmath+1}(s))))]-D^{\alpha}[(v(s )- \psi (s , v(s ),v(v(s))))]$$
$$\leq \aleph_{1}(s ,\sigma_{\jmath}(s),\sigma(\sigma_{\jmath}(s)))+\aleph_{2}(s,\rho_{\jmath}(s),\rho(\rho_{\jmath}(s)))-\aleph_{1}(s , v(s ),v(v(s)))-\aleph_{2}(s , v(s ),v(v(s)))$$
$$\leq 0$$
for the whole $s\in I$ and $\Theta(s_{ 0})=0$. This yields
$$\sigma_{\jmath+1}(s)-\psi (s ,\sigma_{\jmath+1}(s),\sigma(\sigma_{\jmath+1}(s)))\leq v(s )- \psi (s , v(s ),v(v(s)))$$

for every $s\in I$. Since assumption (a1) is valid, we get
$\sigma_{\jmath+1}(s)\leq v(s)$, $\forall s\in I$. Likewise it can
be demonstrated which $v(s)\leq \rho_{\jmath+1}(s )$ on $I$. In
principle, the method of induction, $\sigma_{t}\leq v \leq \rho_{t}$
for every $s\in I$. By taking $t\rightarrow \infty $ limit, we get
$\sigma\leq v \leq\rho$ on $I$. So $(\sigma,\rho)$  they are mixed
type (a) extreme solutions for the equations \eqref{hde} on $I$.,
i.e,

\begin{equation}
\left\{ \begin{array}{c}
     D^{\alpha}[\sigma(s)-\psi(s,\sigma(s),\sigma(\sigma(s))]\leq \aleph_{1}(s, \sigma(s), \sigma(\sigma(s)))
     +\aleph_{1}(s,\rho(s),\rho(\rho(s)))),   s\in
I,\\
        \sigma(s_{0})= v_{0}
        \end{array}\right.\label{hde7}
\end{equation}
and

\begin{equation}
\left\{ \begin{array}{c}
     D^{\alpha}[\rho(s)-\psi(s,\rho(s),\rho(\rho(s))]\geq \aleph_{1}(s, \rho(s), \rho(\rho(s)))+\aleph_{1}(s,\sigma(s),\sigma(\sigma(s))),   s\in
I,\\
        \rho(s_{0})= v_{0}
        \end{array}\right.\label{hde54}
\end{equation}

 the proof is completed .\,\,\,\,\,\,\,\,\,\,\,\,\,\,\,\,\,\,\,\,\,\,\,\,\,\,\,\,\,\,\,\,\,\,\,\,\,\,$\square$


\bigskip \noindent \textbf{Corollary 3.1}\label{re2}\,
Suppose the hypothesis of Theorem 3.4 are fulfilled. Assume that for
$\imath _{1}\geq \imath_{2}$, $\imath _{1},\imath_{2} \in
\overline{\mho}$, then
$$\aleph_{1}(s, \imath_{1}(s), \imath(\imath_{1}(s)))-\aleph_{1}(s,\imath_{2}(s),\imath(\imath_{2}(s)))\leq N_{1}[\imath_{1}(s)-\psi(s,\imath_{1}(s),\imath(\imath_{1}(s)))-(\imath_{2}(s)-\psi(s,\imath_{2}(s),\imath(\imath_{2}(s))), N_{1}> 0,$$
and
$$\aleph_{2}(s, \imath_{1}(s), \imath(\imath_{1}(s)))-\aleph_{2}(s,\imath_{2}(s),\imath(\imath_{2}(s)))\leq N_{2}[\imath_{1}(s)-\psi(s,\imath_{1}(s),\imath(\imath_{1}(s)))-(\imath_{2}(s)-\psi(s,\imath_{2}(s),\imath(\imath_{2}(s))), N_{2}> 0,$$
thus $\sigma(s)= v(s)= \rho(s)$ on $I$ .

\bigskip \noindent \textbf{Proof.}
For $\sigma \leq \rho$ on $I$, it suffices to demonstrate that
$\rho\leq \sigma$ on $I$. Introduce  a function $\Theta \in C
(I,\Re)$
$$\Theta(s )- \psi(s ,\Theta(s),\Theta(\Theta(s)))= (\rho(s)-\psi(s ,\rho(s),\rho(\rho(s))))-(\sigma(s)-\psi(s,\sigma(s),\sigma(\sigma(s)))).$$
Next, $\Theta(s_{0}) =0$ and
$$D^{\alpha}[\Theta(s )- \psi(s ,\Theta(s),\Theta(\Theta(s)))]= D^{\alpha}[(\rho(s)-\psi(s ,\rho(s),\rho(\rho(s))))]-D^{\alpha}[(\sigma(s)-\psi(s,\sigma(s),\sigma(\sigma(s))))]$$
$$= \aleph_{1}(s ,\rho(s),\rho(\rho(s)))- \aleph_{1}(s,\sigma(s),\sigma(\sigma(s)))+ \aleph_{2}(s,\sigma(s),\sigma(\sigma(s)))- \aleph_{2}(s ,\rho(s),\rho(\rho(s)))$$
$$\leq N_{1}[(\rho(s)- \psi(s ,\rho(s),\rho(\rho(s)))-(\sigma(s)-\psi(s,\sigma(s),\sigma(\sigma(s))))]$$

$$+ N_{2}[(\sigma(s)- \sigma(s,\sigma(s),\sigma(\sigma(s))))-(\rho(s)- \psi(s ,\rho(s),\rho(\rho(s))))]$$
$$= (N_{1}+N_{2})[\Theta(s )- \psi(s ,\Theta(s),\Theta(\Theta(s)))].$$
This demonstrates that $\Theta(s )- \psi(s ,\Theta(s),\Theta(\Theta(s)))\leq 0$ on $I$, demonstrating that $\rho \leq \sigma$ on $I$. Therefore $\sigma = \rho = v$ $I$, the proof is completed.\,\,\,\,\,\,\,\,\,\,\,\,\,\,\,\,\,\,\,\,\,\,\,\,\,\,\,\,\,\,\,\,\,\,\,\,\,\,$\square$



\bigskip \noindent \textbf{Theorem 3.5}\label{reR2}\,
Let us suppose that the assumption $(a1)-(a2)$ and $(b2) - (b4)$ achieved. Therefore, for any solution
$v(s)$ of \eqref{hde} with $\sigma_{0} \leq v \leq \rho_{0}$, and we are an iteration $\sigma_{t}, \rho_{t}$ satisfactory for $s\in I$,
\begin{equation}
\left\{ \begin{array}{c}
     \sigma_{0} \leq \sigma_{2} \leq \ldots \leq \sigma_{2t} \leq v \leq \sigma_{2t+1} \leq \ldots \leq \sigma_{3} \leq \sigma_{1},\\
     \rho_{1} \leq \rho_{3} \leq \ldots \leq \rho_{2t+1} \leq v \leq \rho_{2t} \leq \ldots \leq \rho_{2} \leq \rho_{0},
\end{array}\right.\label{d}
\end{equation}
as long as $\sigma_{0} \leq \sigma_{2}$ and $ \rho_{2}\leq\rho_{0}$ on $I$, in which iterating is given by
\begin{equation}
\left\{ \begin{array}{c}
     D^{\alpha}[\sigma_{2t+1} (s)- \psi(s ,\sigma_{2t+1}(s),\sigma(\sigma_{2t+1}(s))]= \aleph_{1}(s,\rho_{t}(s),\rho(\rho_{t}(s)))+\aleph_{2}(s,\sigma_{t}(s),\sigma(\sigma_{t}(s))),\,\,\, s \in I,\\
     \sigma_{2t+1}(s_{0})=v_{o},
\end{array}\right.\label{r}
\end{equation}
and
\begin{equation}
\left\{ \begin{array}{c}
     D^{\alpha}[\rho_{2t+1}(s)- \psi(s ,\rho_{2t+1}(s),\rho(\rho_{2t+1}(s))]= \aleph_{1}(s,\sigma_{t}(s),\sigma(\sigma_{t}(s)))+\aleph_{2}(s,\rho_{t}(s),\rho\rho_{t}(s))),\,\,\, s \in I,\\
     \rho_{2t+1}(s_{0})=v_{o},
\end{array}\right.\label{1r}
\end{equation}

of $t\in N$. Furthermore, the monotonous sequences $\{\sigma_{2t}\},\{\sigma_{2t+1}\}, \{\rho_{2t}\},\{\rho_{2t+1}\}$ converge uniformly to $\sigma, \rho, \sigma^{\diamond}, \rho^{\diamond}$, respectively, and fulfilling this assumptions:

$$(1) \,\, D^{\alpha}[\sigma(s)- \psi(s ,\sigma(s),\sigma(\sigma(s)))]= \aleph_{1}(s,\rho(s),\rho(\rho(s)))+\aleph_{2}(s,\sigma(s),\sigma(\sigma(s)))$$
$$(2) \,\,\, D^{\alpha}[\rho(s)- \psi(s ,\rho(s),\rho(\rho(s))]= \aleph_{1}(s,\sigma_{t}(s),\sigma(\sigma(s)))+\aleph_{2}(s,\rho(s),\rho\rho(s)))$$
$$(3)\,\,\, D^{\alpha}[\sigma^{\diamond}(s)- \psi(s ,\sigma^{\diamond}(s),\sigma(\sigma^{\diamond}(s)))]= \aleph_{1}(s,\rho^{\diamond}(s),\rho(\rho^{\diamond}(s)))+\aleph_{2}(s,\sigma^{\diamond}(s),\sigma(\sigma^{\diamond}(s)))$$
$$(4) \,\,\, D^{\alpha}[\rho^{\diamond}(s)- \psi(s ,\rho^{\diamond}(s),\rho(\rho^{\diamond}(s))]=\aleph_{1}(s,\sigma^{\diamond}(s),\sigma(\sigma^{\diamond}(s)))+\aleph_{2}(s,\rho^{\diamond}(s),\rho\rho^{\diamond}(s)))$$
$for s \in I and \sigma \leq v \leq \rho, \sigma^{\diamond}\leq v \leq \rho^{\diamond}, s \in I , \sigma(0) =\sigma(0) = \sigma^{\diamond}(0) = \rho^{\diamond}(0) = v0.$\\

\bigskip \noindent \textbf{Proof.} By the assumptions of the
theorem, we suppose that $\sigma_{0}\leq \sigma_{2}$ and
$\rho_{2}\leq\rho_{0} $, on $I$. We demonstrate that
\begin{equation}
\left\{\begin{array}{c}
    \sigma_{0}\leq \sigma_{2} \leq v \leq \sigma_{3}\leq \sigma_{1},\\
    \rho_{1}\leq\rho_{3} \leq v \leq \rho_{2}\leq\rho_{0}
\end{array}\right.\label{2r}
\end{equation}
on $I$. Set
$$\Theta(s )- \psi(s ,\Theta(s ),\Theta(\Theta(s )))= (v(s )- \psi(s,v(s ),v(v(s ))))-(\sigma_{1}(s)- \psi (s ,\sigma_{1}(s),\sigma(\sigma_{1}(s)))).$$
utilization that $\sigma_{ 0 }\leq v \leq \rho_{0}$ on $I$, as $v$ is any solution of
 \eqref{hde} and the monotonous the nature of functions $\aleph_{1}$ and
 $\aleph_{2}$, this yields

$$D^{\alpha}[\Theta(s )- \psi(s ,\Theta(s ),\Theta(\Theta(s )))]= D^{\alpha}[(v(s )-
\psi(s,v(s ),v(v(s ))))]-D^{\alpha}[(\sigma_{1}(s)- \psi (s ,\sigma_{1}(s),\sigma(\sigma_{1}(s))))]$$
$$=\aleph_{1}(s,v(s ),v(v(s )))))+\aleph_{2}(s,v(s ),v(v(s ))))-\aleph_{1}(s,\rho_{0}(s),\rho(\rho_{0}(s)))
-\aleph_{2}(s ,\sigma_{0}(s),\sigma(\sigma_{0}(s))))$$
$$\leq0$$

for every $s \in I$ and $\Theta(s_{0}) = 0$. Thus, we reached the conclusion
$$v(s )- \psi(s,v(s ),v(v(s )))\leq \sigma_{1}(s)- \psi (s ,\sigma_{1}(s),\sigma(\sigma_{1}(s)))$$or
$$v(s )\leq \sigma_{1}(s)$$
for every $s \in I.$ In the same way,  we can show that
$\sigma_{3}\leq\sigma_{1}, \rho_{1} \leq v $ and $\sigma_{2} \leq
v$, taking into account differences
$$\Theta(s)-\psi(s,\Theta(s),\Theta(\Theta(s)))= (\sigma_{3}(s)-\psi(s,\sigma_{3}(s),\sigma( \sigma_{3}(s))))-
( \sigma_{1}(s)-\psi(s ,\sigma_{1}(s),\sigma(\sigma_{1}(s)))),$$
$$\Theta(s)-\psi(s,\Theta(s),\Theta(\Theta(s)))= (\rho_{1}(s)-\psi(s,\rho_{1}(s),\rho(\rho_{1}(s))))-( v(s)-\psi(s ,v(s),v(v(s))))$$
and
$$\Theta(s)-\psi(s,\Theta(s),\Theta(\Theta(s)))= (\sigma_{2}(s)-\psi(s,\sigma_{2}(s),\sigma( \sigma_{2}(s))))-( v(s)-\psi(s ,v(s),v(v(s))))$$

respectively. At each of these cases, we get $\Theta(s)
-\psi(s,\Theta(s),\Theta(\Theta(s))\leq 0$, for all $s \in I$ and
representation \eqref{2r} is established. This completed prove.


\bigskip \noindent \textbf{ Competing Interests}  The authors declare that
   they have no competing interests.

 \bigskip \noindent \textbf{ Authors$'$ contributions}
All the authors jointly worked on deriving the results and approved
the final manuscript.


\bigskip
\small

\begin{thebibliography}{10}
%






\bibitem{loverro2004fractional}  Loverro A. (2004). Fractional Calculus: History, Definitions and Applications for the Engineer

\bibitem{podlubny1998fractional} Podlubny, I. (1998). Fractional differential equations: an introduction to fractional derivatives,
fractional differential equations, to methods of their solution and  some of their applications (Vol. 198). Academic press.

\bibitem{millerintroduction}Miller, K. S., and Ross, B. An Introduction to the Fractional Calculus and Differential Equations. 1993.


\bibitem{kulish2002application} Kulish, V. V., and Lage, J. L. (2002). Application of fractional calculus to fluid mechanics. Journal of Fluids Engineering, 124(3), 803-806.

\bibitem{havlin1995fractals} Havlin, S., Buldyrev, S. V., Goldberger, A. L., Mantegna, R. N., Ossadnik, S. M., Peng, C. K., ... and Stanley, H. E. (1995). Fractals in biology and medicine. Chaos, Solitons and Fractals, 6, 171-201.

\bibitem{nonnenmacher2013fractals} Nonnenmacher, T. F., Losa, G. A., and Weibel, E. R. (Eds.). (2013). Fractals in biology and medicine. Birkh$\ddot{ä}$user.


\bibitem{smith2007bacterial} Smith, H. L. (2007). Bacterial Growth. Retrieved on, 09-15.


\bibitem{dhage2012basic} Dhage, B. C. (2012).
 Basic results in the theory of hybrid differential equations with linear perturbations os second type.
 Tamkang Journal of Mathematics, 44(2), 171-186.

\bibitem{lu2013theory} Lu, H., Sun, S., Yang, D., and Teng, H.(2013). Theory of fractional hybrid differential equations with linear perturbations of second type. Boundary Value Problems, 2013(1), 1-16.


 \bibitem{dhage2014approximation}  Dhage, B. C. (2014). Approximation methods in the theory of
  hybrid differential equations with linear perturbations of second type. Tamkang Journal of Mathematics, 45(1), 39-61.

  \bibitem{dhage2013basic}  Dhage, B. C.,an Jadhav, N. S. (2013). Basic Results in the Theory of Hybrid Differential Equations
  with Linear Perturbations of Second Type. Tamkang Journal of Mathematics, 44(2), 171-186.

\bibitem{Dhage2010414}  Dhage, B. C.,an Lakshmikantham, V. (2010). Basic results on hybrid differential
equations. Nonlinear Analysis: Hybrid Systems, 4(3), 414-424.

\bibitem{ibrahim2012} Ibrahim R.W. (2012). Existence of deviating fractional differential
equation, CUBO A Mathematical Journal, 14 (03), 127-140.

 \bibitem{ibrahim2015existence} Ibrahim, R. W., Kili\c{c}man, A., and Damag, F. H. (2015).
   Existence and uniqueness for a class of iterative fractional differential equations. Advances in Difference Equations, 2015(1), 1-13.


\bibitem{dhage2004fixed} Dhage, B., and O'Regan, D. (2004). A fixed point theorem in Banach algebras with applications to functional integral equations.
 Functional Differential Equations, 7(3-4), p-259.

\bibitem{Granas1991R} Granas A., Guenther R. B., Lee J.W. (1991) Some general existence principles for Caratheodory theory
  of nonlinear differential equations, J. Math. Pures. Appl., 70 , 153-196.


\end{thebibliography}
\end{document}